\documentclass[12pt]{article}%
\usepackage{graphicx}
\usepackage[intlimits]{amsmath}
\usepackage{latexsym}
\usepackage{amsfonts}
\usepackage{amssymb}%
\makeatletter
\newfont{\footsc}{cmcsc10 at 8truept}
\newfont{\footbf}{cmbx10 at 8truept}
\newfont{\footrm}{cmr10 at 10truept}
\newtheorem{theorem}{Theorem}

\newtheorem{lemma}[theorem]{Lemma}

\newtheorem{problem}[theorem]{Problem}

\newenvironment{proof}[1][Proof]{\noindent{\textbf {#1}  }}  {\hfill$\Box$\bigskip}

\begin{document}

\title{Chromatic number and spectral radius}
\author{Vladimir Nikiforov\\{\small Department of Mathematical Sciences, University of Memphis, }\\{\small Memphis TN 38152, USA; e-mail: vnikifrv@memphis.edu}}
\maketitle

\begin{abstract}
Write $\mu\left(  A\right)  =\mu_{1}\left(  A\right)  \geq\cdots\geq\mu_{\min
}\left(  A\right)  $ for the eigenvalues of a Hermitian matrix $A$. Our main
result is: let $A$ be a Hermitian matrix partitioned into $r\times r$ blocks
so that all diagonal blocks are zero. Then for every real diagonal matrix $B$
of the same size as $A,$
\[
\mu\left(  B-A\right)  \geq\mu\left(  B+\frac{1}{r-1}A\right)  .
\]

Let $G$ be a nonempty graph, $\chi\left(  G\right)  $ be its chromatic number,
$A$ be its adjacency matrix, and $L$ be its Laplacian. The above inequality
implies the well-known result of A.J. Hoffman%
\[
\chi\left(  G\right)  \geq1+\frac{\mu\left(  A\right)  }{-\mu_{\min}\left(
A\right)  },
\]
and also,
\[
\chi\left(  G\right)  \geq1+\frac{\mu\left(  A\right)  }{\mu\left(  L\right)
-\mu\left(  A\right)  }.
\]
Equality holds in the latter inequality if and only if every two color classes
of $G\ $induce a $\left\vert \mu_{\min}\left(  A\right)  \right\vert $-regular
subgraph.\bigskip



\end{abstract}

Write $\mu\left(  A\right)  =\mu_{1}\left(  A\right)  \geq\cdots\geq\mu_{\min
}\left(  A\right)  $ for the eigenvalues of a Hermitian matrix $A$. Given a
graph $G,$ let $\chi\left(  G\right)  $ be its chromatic number, $A\left(
G\right)  $ be its adjacency matrix, and $D\left(  G\right)  $ be the diagonal
matrix of its degree sequence; set $L\left(  G\right)  =D\left(  G\right)
-A\left(  G\right)  .$

Letting $G$ be a nonempty graph with $L\left(  G\right)  =L$ and $A\left(
G\right)  =A,$ we prove that%
\begin{equation}
\chi\left(  G\right)  \geq1+\frac{\mu\left(  A\right)  }{\mu\left(  L\right)
-\mu\left(  A\right)  }, \label{mainin}%
\end{equation}
complementing the well-known inequality of A.J. Hoffman \cite{Hof70}
\begin{equation}
\chi\left(  G\right)  \geq1+\frac{\mu\left(  A\right)  }{-\mu_{\min}\left(
A\right)  }. \label{hofin}%
\end{equation}
Equality holds in (\ref{mainin}) if and only if every two color classes of $G$
induce a $\left\vert \mu_{\min}\left(  A\right)  \right\vert $-regular subgraph.

We deduce inequalities (\ref{mainin}) and (\ref{hofin}) from a theorem of its
own interest.

\begin{theorem}
\label{th1} Let $A$ be a Hermitian matrix partitioned into $r\times r$ blocks
so that all diagonal blocks are zero. Then for every real diagonal matrix $B$
of the same size as $A,$
\begin{equation}
\mu\left(  B-A\right)  \geq\mu\left(  B+\frac{1}{r-1}A\right)  .
\label{mbound}%
\end{equation}

\end{theorem}

\begin{proof}
[\textbf{Proof of Theorem \ref{th1}}]Write $n$ for the size of $A,$ let
$\left[  n\right]  =\cup_{i=1}^{r}N_{i}$ be the partition of its index set,
and let $b_{1},\ldots,b_{n}$ be the diagonal entries of $B.$ Set $L=B-A,$
$K=\left(  r-1\right)  B+A,$ and select a unit eigenvector $\mathbf{x}=\left(
x_{1},\ldots,x_{n}\right)  $ to $\mu\left(  K\right)  .$ Our proof strategy is
simple: using $\mathbf{x},$ we define specific $n$-vectors $\mathbf{y}%
_{1},\ldots,\mathbf{y}_{r}$ and show that
\[
r\left(  r-1\right)  \mu\left(  L\right)  \geq\mu\left(  L\right)
%TCIMACRO{\tsum \limits_{i\in\left[  r\right]  }}%
%BeginExpansion
{\textstyle\sum\limits_{i\in\left[  r\right]  }}
%EndExpansion
\left\Vert \mathbf{y}_{i}\right\Vert ^{2}\geq%
%TCIMACRO{\tsum \limits_{i\in\left[  r\right]  }}%
%BeginExpansion
{\textstyle\sum\limits_{i\in\left[  r\right]  }}
%EndExpansion
\left\langle L\mathbf{y}_{i},\mathbf{y}_{i}\right\rangle \geq r\left\langle
K\mathbf{x},\mathbf{x}\right\rangle =r\mu\left(  K\right)  .
\]
For $i=1,\ldots,r$ define $\mathbf{y}_{i}=\left(  y_{i1},\ldots,y_{in}\right)
$ as%
\[
y_{ij}=\left\{
\begin{array}
[c]{ll}%
-\left(  r-1\right)  x_{j} & \text{if }j\in N_{i}\\
x_{j} & \text{if }j\in\left[  n\right]  \backslash N_{i}.
\end{array}
\right.
\]
The Rayleigh principle implies that
\begin{equation}
\mu\left(  L\right)
%TCIMACRO{\tsum \limits_{i\in\left[  r\right]  }}%
%BeginExpansion
{\textstyle\sum\limits_{i\in\left[  r\right]  }}
%EndExpansion
\left\Vert \mathbf{y}_{i}\right\Vert ^{2}\geq%
%TCIMACRO{\tsum \limits_{i\in\left[  r\right]  }}%
%BeginExpansion
{\textstyle\sum\limits_{i\in\left[  r\right]  }}
%EndExpansion
\mu\left(  L\right)  \left\Vert \mathbf{y}_{i}\right\Vert ^{2}\geq%
%TCIMACRO{\tsum \limits_{i\in\left[  r\right]  }}%
%BeginExpansion
{\textstyle\sum\limits_{i\in\left[  r\right]  }}
%EndExpansion
\left\langle L\mathbf{y}_{i},\mathbf{y}_{i}\right\rangle . \label{in1}%
\end{equation}
Noting that
\[
\left\Vert \mathbf{y}_{i}\right\Vert ^{2}=%
%TCIMACRO{\tsum \limits_{j\in\left[  n\right]  \backslash N_{i}}}%
%BeginExpansion
{\textstyle\sum\limits_{j\in\left[  n\right]  \backslash N_{i}}}
%EndExpansion
\left\vert x_{j}\right\vert ^{2}+\left(  r-1\right)  ^{2}%
%TCIMACRO{\tsum \limits_{j\in N_{i}}}%
%BeginExpansion
{\textstyle\sum\limits_{j\in N_{i}}}
%EndExpansion
\left\vert x_{j}\right\vert ^{2}=1+r\left(  r-2\right)
%TCIMACRO{\tsum \limits_{j\in N_{i}}}%
%BeginExpansion
{\textstyle\sum\limits_{j\in N_{i}}}
%EndExpansion
\left\vert x_{j}\right\vert ^{2},
\]
we obtain,%
\begin{equation}%
%TCIMACRO{\tsum \limits_{i\in\left[  r\right]  }}%
%BeginExpansion
{\textstyle\sum\limits_{i\in\left[  r\right]  }}
%EndExpansion
\left\Vert \mathbf{y}_{i}\right\Vert ^{2}=r+r\left(  r-2\right)
%TCIMACRO{\tsum \limits_{i\in\left[  r\right]  }}%
%BeginExpansion
{\textstyle\sum\limits_{i\in\left[  r\right]  }}
%EndExpansion%
%TCIMACRO{\tsum \limits_{j\in N_{i}}}%
%BeginExpansion
{\textstyle\sum\limits_{j\in N_{i}}}
%EndExpansion
\left\vert x_{j}\right\vert ^{2}=r+r\left(  r-2\right)  =r\left(  r-1\right)
. \label{in2}%
\end{equation}
On the other hand, we have
\[
\left\langle L\mathbf{y}_{i},\mathbf{y}_{i}\right\rangle =%
%TCIMACRO{\tsum \limits_{j\in\left[  n\right]  }}%
%BeginExpansion
{\textstyle\sum\limits_{j\in\left[  n\right]  }}
%EndExpansion
b_{j}\left\vert y_{ij}\right\vert ^{2}-%
%TCIMACRO{\tsum \limits_{j,k\in\left[  n\right]  }}%
%BeginExpansion
{\textstyle\sum\limits_{j,k\in\left[  n\right]  }}
%EndExpansion
a_{jk}y_{ik}\overline{y_{ij}}.
\]
For every $i\in\left[  n\right]  ,$ we see that
\[%
%TCIMACRO{\tsum \limits_{j\in\left[  n\right]  }}%
%BeginExpansion
{\textstyle\sum\limits_{j\in\left[  n\right]  }}
%EndExpansion
b_{j}\left\vert y_{ij}\right\vert ^{2}=%
%TCIMACRO{\tsum \limits_{j\in\left[  n\right]  }}%
%BeginExpansion
{\textstyle\sum\limits_{j\in\left[  n\right]  }}
%EndExpansion
b_{j}\left\vert x_{j}\right\vert ^{2}+r\left(  r-2\right)
%TCIMACRO{\tsum \limits_{j\in N_{i}}}%
%BeginExpansion
{\textstyle\sum\limits_{j\in N_{i}}}
%EndExpansion
b_{j}\left\vert x_{j}\right\vert ^{2},
\]
and, likewise,
\[%
%TCIMACRO{\tsum \limits_{j,k\in\left[  n\right]  }}%
%BeginExpansion
{\textstyle\sum\limits_{j,k\in\left[  n\right]  }}
%EndExpansion
a_{jk}y_{ik}\overline{y_{ij}}=%
%TCIMACRO{\tsum \limits_{j,k\in\left[  n\right]  }}%
%BeginExpansion
{\textstyle\sum\limits_{j,k\in\left[  n\right]  }}
%EndExpansion
a_{jk}x_{k}\overline{x_{j}}-r%
%TCIMACRO{\tsum \limits_{j\in N_{i},k\in\left[  n\right]  }}%
%BeginExpansion
{\textstyle\sum\limits_{j\in N_{i},k\in\left[  n\right]  }}
%EndExpansion
a_{jk}x_{k}\overline{x_{j}}-r%
%TCIMACRO{\tsum \limits_{k\in N_{i},j\in\left[  n\right]  }}%
%BeginExpansion
{\textstyle\sum\limits_{k\in N_{i},j\in\left[  n\right]  }}
%EndExpansion
a_{jk}x_{k}\overline{x_{j}}.
\]
Summing these for all $i\in\left[  r\right]  ,$ we find that
\begin{align*}%
%TCIMACRO{\tsum \limits_{i\in\left[  r\right]  }}%
%BeginExpansion
{\textstyle\sum\limits_{i\in\left[  r\right]  }}
%EndExpansion
\left\langle L\mathbf{y}_{i},\mathbf{y}_{i}\right\rangle  &  =%
%TCIMACRO{\tsum \limits_{i\in\left[  r\right]  ,j\in\left[  n\right]  }}%
%BeginExpansion
{\textstyle\sum\limits_{i\in\left[  r\right]  ,j\in\left[  n\right]  }}
%EndExpansion
b_{j}\left\vert x_{j}\right\vert ^{2}+r\left(  r-2\right)
%TCIMACRO{\tsum \limits_{i\in\left[  r\right]  ,j\in N_{i}}}%
%BeginExpansion
{\textstyle\sum\limits_{i\in\left[  r\right]  ,j\in N_{i}}}
%EndExpansion
b_{j}\left\vert x_{j}\right\vert ^{2}\\
&  -r%
%TCIMACRO{\tsum \limits_{j,k\in\left[  n\right]  }}%
%BeginExpansion
{\textstyle\sum\limits_{j,k\in\left[  n\right]  }}
%EndExpansion
a_{jk}x_{k}\overline{x_{j}}+r%
%TCIMACRO{\tsum \limits_{i\in\left[  r\right]  }}%
%BeginExpansion
{\textstyle\sum\limits_{i\in\left[  r\right]  }}
%EndExpansion
\left(
%TCIMACRO{\tsum \limits_{j\in N_{i},k\in\left[  n\right]  }}%
%BeginExpansion
{\textstyle\sum\limits_{j\in N_{i},k\in\left[  n\right]  }}
%EndExpansion
a_{jk}x_{k}\overline{x_{j}}+r%
%TCIMACRO{\tsum \limits_{k\in N_{i},j\in\left[  n\right]  }}%
%BeginExpansion
{\textstyle\sum\limits_{k\in N_{i},j\in\left[  n\right]  }}
%EndExpansion
a_{jk}x_{k}\overline{x_{j}}\right) \\
&  =r\left(  r-1\right)
%TCIMACRO{\tsum \limits_{j\in\left[  n\right]  }}%
%BeginExpansion
{\textstyle\sum\limits_{j\in\left[  n\right]  }}
%EndExpansion
b_{j}\left\vert x_{j}\right\vert ^{2}-r%
%TCIMACRO{\tsum \limits_{j,k\in\left[  n\right]  }}%
%BeginExpansion
{\textstyle\sum\limits_{j,k\in\left[  n\right]  }}
%EndExpansion
a_{jk}x_{k}\overline{x_{j}}+2r%
%TCIMACRO{\tsum \limits_{j,k\in\left[  n\right]  }}%
%BeginExpansion
{\textstyle\sum\limits_{j,k\in\left[  n\right]  }}
%EndExpansion
a_{jk}x_{k}\overline{x_{j}}\\
&  =r\left(  r-1\right)
%TCIMACRO{\tsum \limits_{j\in\left[  n\right]  }}%
%BeginExpansion
{\textstyle\sum\limits_{j\in\left[  n\right]  }}
%EndExpansion
b_{j}\left\vert x_{j}\right\vert ^{2}+r%
%TCIMACRO{\tsum \limits_{j,k\in\left[  n\right]  }}%
%BeginExpansion
{\textstyle\sum\limits_{j,k\in\left[  n\right]  }}
%EndExpansion
a_{jk}x_{k}\overline{x_{j}}=r\left\langle K\mathbf{x},\mathbf{x}\right\rangle
=r\mu\left(  K\right)  .
\end{align*}
Hence, in view of (\ref{in1}) and (\ref{in2}), we obtain $\left(  r-1\right)
\mu\left(  B-A\right)  \geq\mu\left(  K\right)  ,$ completing the proof.
\end{proof}

\bigskip

\begin{lemma}
\label{le1}Let $A$ be an irreducible nonnegative symmetric matrix and $R$ be
the diagonal matrix of its rowsums. Then
\begin{equation}
\mu\left(  R+\frac{1}{r-1}A\right)  \geq\frac{r}{r-1}\mu\left(  A\right)
\label{ineq}%
\end{equation}
with equality holding if and only if all rowsums of $A$ are equal.
\end{lemma}

\begin{proof}
Let $A=\left(  a_{ij}\right)  $ and $n$ be its size. Note first that for any
vector $\mathbf{x}=\left(  x_{1},\ldots,x_{n}\right)  $
\[
\left\langle \left(  R-A\right)  \mathbf{x},\mathbf{x}\right\rangle =%
%TCIMACRO{\tsum \limits_{1\leq i<j\leq n}}%
%BeginExpansion
{\textstyle\sum\limits_{1\leq i<j\leq n}}
%EndExpansion
a_{ij}\left(  x_{i}-x_{j}\right)  ^{2}\geq0.
\]
Hence, $R-A$ is positive semidefinite; since $A$ is irreducible, if
$\left\langle \left(  R-A\right)  \mathbf{x},\mathbf{x}\right\rangle =0,$ then
all entries of $\mathbf{x}$ are equal.

Let $\mathbf{x}=\left(  x_{1},\ldots,x_{n}\right)  $ be an eigenvector to
$\mu=\mu\left(  D+\frac{1}{r-1}A\right)  .$ We have
\begin{align}
\mu &  =\left\langle \left(  R+\frac{1}{r-1}A\right)  \mathbf{x}%
,\mathbf{x}\right\rangle =%
%TCIMACRO{\tsum \limits_{i=1}^{n}}%
%BeginExpansion
{\textstyle\sum\limits_{i=1}^{n}}
%EndExpansion
x_{i}^{2}%
%TCIMACRO{\tsum \limits_{j=1}^{n}}%
%BeginExpansion
{\textstyle\sum\limits_{j=1}^{n}}
%EndExpansion
a_{ij}+\frac{1}{r-1}%
%TCIMACRO{\tsum \limits_{i=1}^{n}}%
%BeginExpansion
{\textstyle\sum\limits_{i=1}^{n}}
%EndExpansion%
%TCIMACRO{\tsum \limits_{j=1}^{n}}%
%BeginExpansion
{\textstyle\sum\limits_{j=1}^{n}}
%EndExpansion
a_{ij}x_{i}x_{j}\nonumber\\
&  =%
%TCIMACRO{\tsum \limits_{1\leq i<j\leq n}}%
%BeginExpansion
{\textstyle\sum\limits_{1\leq i<j\leq n}}
%EndExpansion
a_{ij}\left(  x_{i}-x_{j}\right)  ^{2}+\frac{r}{r-1}%
%TCIMACRO{\tsum \limits_{i=1}^{n}}%
%BeginExpansion
{\textstyle\sum\limits_{i=1}^{n}}
%EndExpansion%
%TCIMACRO{\tsum \limits_{j=1}^{n}}%
%BeginExpansion
{\textstyle\sum\limits_{j=1}^{n}}
%EndExpansion
a_{ij}x_{i}x_{j}\nonumber\\
&  =\left\langle \left(  R-A\right)  \mathbf{x},\mathbf{x}\right\rangle
+\frac{r}{r-1}%
%TCIMACRO{\tsum \limits_{i=1}^{n}}%
%BeginExpansion
{\textstyle\sum\limits_{i=1}^{n}}
%EndExpansion%
%TCIMACRO{\tsum \limits_{j=1}^{n}}%
%BeginExpansion
{\textstyle\sum\limits_{j=1}^{n}}
%EndExpansion
a_{ij}x_{i}x_{j}\geq\frac{r}{r-1}\mu\left(  A\right)  , \label{in_3}%
\end{align}
proving (\ref{ineq}).

Let now equality holds in (\ref{ineq}). Then equality holds in (\ref{in_3}),
and so $\left\langle \left(  R-A\right)  \mathbf{x},\mathbf{x}\right\rangle
=0$ and $\mathbf{x}$ is an eigenvector of $A$ to $\mu\left(  A\right)  .$
Therefore $x_{1}=\cdots=x_{n}$ and the rowsums of $A$ are equal.

If the rowsums of $A$ are equal, the vector $\mathbf{j}=\left(  1,\ldots
,1\right)  $ is an eigenvector of $A$ to $\mu\left(  A\right)  $ and of $R$ to
$\mu\left(  R\right)  ;$ therefore $\mathbf{j}$ is an eigenvector of
$R+\frac{1}{r-1}A$ to $\mu$, and so equality holds in (\ref{ineq}), completing
the proof.
\end{proof}

\bigskip

\begin{proof}
[\textbf{Proof of (\ref{mainin}) and (\ref{hofin})}]Let $G$ be a graph with
chromatic number $\chi=r.$ Coloring the vertices of $G$ into $r$ colors
defines a partition of its adjacency matrix $A=A\left(  G\right)  $ with zero
diagonal blocks. Letting $B$ be the zero matrix, Theorem \ref{th1} implies
inequality (\ref{hofin}).

Letting now $B=D=D\left(  G\right)  ,$ Lemma \ref{le1} implies that
\[
\mu\left(  D+\frac{1}{r-1}A\right)  \geq\frac{r}{r-1}\mu\left(  A\right)  ,
\]
and inequality (\ref{mainin}) follows.

The following argument for equality in (\ref{mainin}) was kindly suggested by
the referee. If equality holds in (\ref{mainin}), by Lemma \ref{le1}, $G$ is
regular; hence equality holds also in (\ref{hofin}). Setting $\mu\left(
G\right)  =k,$ $\left\vert \mu_{\min}\left(  G\right)  \right\vert =\tau$ and
writing $\alpha\left(  G\right)  $ for the independence number of $G,$ let us
recall Hoffman's bound on $\alpha\left(  G\right)  $: for every $k$-regular
graph $G,$%
\begin{equation}
\alpha\left(  G\right)  \leq\frac{n\tau}{k+\tau}. \label{hofup}%
\end{equation}
On the other hand, we have
\[
\alpha\left(  G\right)  \geq\frac{n}{\chi\left(  G\right)  }=\frac{n}%
{1+k/\tau}=\frac{n\tau}{k+\tau}.
\]
and thus, equality holds in (\ref{hofup}). It is known (see, e.g.,
\cite{GoRo01}, Lemma 9.6.2) that this is only possible if $\chi\left(
G\right)  =n/\alpha\left(  G\right)  $ and every two color classes of $G$
induce a $\tau$-regular bipartite subgraph.
\end{proof}

\bigskip

\textbf{Concluding remarks}

For the complete graph of order $n$ without an edge, inequality (\ref{mainin})
gives $\chi=n-1,$ while (\ref{hofin}) gives only $\chi\geq n/2+2.$ By
contrast, for a sufficiently large wheel $W_{1,n}$, i.e., a vertex joined to
all vertices of a cycle of length $n$, we see that (\ref{mainin}) gives
$\chi\geq2,$ while (\ref{hofin}) gives $\chi\geq3.$

A natural question is to determine when equality holds in (\ref{mbound}). A
particular answer, building upon \cite{Mer94}, can be found in \cite{SHW02}:
if $G$ is a connected graph, then $\mu\left(  D\left(  G\right)  -A\left(
G\right)  \right)  =\mu\left(  D\left(  G\right)  +A\left(  G\right)  \right)
$ if and only if $G$ is bipartite.

\begin{problem}
Determine when equality holds in (\ref{mbound}).
\end{problem}

Finally, any lower bound on $\mu\left(  A\left(  G\right)  \right)  ,$
together with (\ref{mainin}), gives a lower bound on $\mu\left(  L\left(
G\right)  \right)  .$ This approach helps deduce some inequalities for
bipartite graphs given in \cite{HoZh05} and \cite{YLT04}.\bigskip

\textbf{Acknowledgement. }Thanks to Peter Rowlinson, Sebi Cioab\u{a} and Cecil
Rousseau for useful suggestions. The author is most indebted to the referee
for the exceptionally thorough, helpful and kind report.

\end{document}